\documentstyle{amlts}
\begin{document}
\annalsline{159}{2004}
\received{November 14, 2002}
\startingpage{447}
\def\bye{\end{document}}
 \font\tenrm=cmr10
\def\ritem#1{\item[{\rm #1}]}
\input amssym.def
\input amssym.tex
\def\eqref#1{(\ref{#1})}
\def\varprojlim{\mathop{\displaystyle\lim_{\longleftarrow}}}
\def\joinrel{\mathrel{\mkern-4mu}}
\def\relbar{\mathrel{\smash-}}
\def\lrar{\relbar\joinrel\relbar\joinrel\relbar\joinrel\relbar\joinrel\rightarrow}
\input boxedeps.tex 
\SetepsfEPSFSpecial 
\HideDisplacementBoxes
\def\figin#1#2{
$$
 {\BoxedEPSF{#1.eps scaled
#2}%
}%
$$
\noindent}
\catcode`\@=11
\font\twelvemsb=msbm10 scaled 1100
\font\tenmsb=msbm10
\font\ninemsb=msbm10 scaled 800
\newfam\msbfam
\textfont\msbfam=\twelvemsb  \scriptfont\msbfam=\ninemsb
  \scriptscriptfont\msbfam=\ninemsb
\def\msb@{\hexnumber@\msbfam}
\def\Bbb{\relax\ifmmode\let\next\Bbb@\else
 \def\next{\errmessage{Use \string\Bbb\space only in math
mode}}\fi\next}
\def\Bbb@#1{{\Bbb@@{#1}}}
\def\Bbb@@#1{\fam\msbfam#1}
\catcode`\@=12

 \catcode`\@=11
\font\twelveeuf=eufm10 scaled 1100
\font\teneuf=eufm10
\font\nineeuf=eufm7 scaled 1100
\newfam\euffam
\textfont\euffam=\twelveeuf  \scriptfont\euffam=\teneuf
  \scriptscriptfont\euffam=\nineeuf
\def\euf@{\hexnumber@\euffam}
\def\frak{\relax\ifmmode\let\next\frak@\else
 \def\next{\errmessage{Use \string\frak\space only in math
mode}}\fi\next}
\def\frak@#1{{\frak@@{#1}}}
\def\frak@@#1{\fam\euffam#1}
\catcode`\@=12
 \newfont{\cyrr}{wncyr10}
\font\mathscr=eusm10 scaled 1100
\font\smathsc=eusm8
\def\Sh{\mbox{\cyrr Sh}}

\def\Q{{\bf Q}}
\def\Z{{\bf Z}}
\def\F{{\bf F}}
\def\R{{\bf R}}
\def\CC{{\bf C}}
\def\v{{\bf v}}

\def\O{{\cal O}}
\def\K{\hbox{\mathscr K}}
\def\ssK{\hbox{\smathsc K}}
\def\A{{\cal A}}
\def\X{{\cal X}}
\def\E{{\cal E}}
\def\C{{\cal C}}
\def\U{{\cal U}}
\def\V{{\cal V}}
\def\L{{\cal L}}

\def\AA{{\frak A}}
\def\p{{\frak p}}

\def\ld{{\cal h}}
\def\rd{{\cal i}}

\def\m{{\frak m}}

\def\bmu{{\mathbold{\mu}}}

\def\Kp{K_p}
\def\Op{\O_p}
\def\Zp{\Z_p}
\def\Fp{\F_p}
\def\Qp{\Q_p}
\def\Opx{\Op^\times}
\def\Qinf{\Q_\infty}
\def\Kinf{K_\infty}
\def\Kinfp{K_{\infty,p}}
\def\Qn{\Q_n}
\def\Qnp{\Q_{n,p}}
\def\Qinfp{\Q_{\infty,p}}
\def\KQinf{{K_{{\rm cyc}}}}
\def\KQinfp{{K_{{\rm cyc},p}}}

\def\too{\longrightarrow}
\def\onto{\twoheadrightarrow}

\def\tw#1{#1(\rho^{-1})}
\def\redgen#1#2{{#2^{\rho}_{#1}}}
\def\rc#1{{#1^{\rho}_{\Kinf}}}
\def\red#1{{#1^{\rho}_{\KQinf}}}
\def\redd#1{\widetilde{#1}}

\def\cp{{\rm char}}
\def\Aut{{\rm Aut}}
\def\Hom{{\rm Hom}}
\def\End{{\rm End}}
\def\Gal{{\rm Gal}}
\def\rank{{\rm rank}}
\def\ord{{\rm ord}}
\def\Sel{{\rm Sel}_{p}}
\def\Selgp{{\rm Sel}_{\p}}
\def\Tr{{\mathbf{Tr}}}
\def\image{{\rm image}}

\def\Iw#1{\Lambda(#1)}
\def\IwO#1{\Lambda_{\O}(#1)}
\def\Iwq{\Lambda}
\def\IwOq{\Lambda_{\O}}

\def\Efg{\hat{E}}
\def\logE{\lambda_E}
\def\map#1{\;\stackrel{#1}{\longrightarrow}\;}
\def\isom{{\map{\,\sim\,}}}
\def\gaminus{\gamma_\ast}
\def\one{\chi_0}

\title{The main 
conjecture for CM elliptic\\ curves at supersingular primes} 
\shorttitle{Main conjecture for supersingular primes} 

  \acknowledgements{The first author was supported by an NSF
Postdoctoral Fellowship.  The second author was supported by NSF grant DMS-0140378.\hfill\break
\hglue23pt 2000 {\it Mathematics Subject Classification}. Primary 11G05, 11R23; Secondary 11G40.} 
 \twoauthors{Robert Pollack}{Karl Rubin}
 \institutions{University of Chicago, Chicago IL\\
{\eightpoint {\it E-mail address\/}: pollack@math.uchicago.edu
}
\\
\vglue6pt
Stanford University, Stanford CA\\
{\eightpoint {\it E-mail address\/}: rubin@math.stanford.edu}} 
\vfil
\centerline{\bf Abstract}
\vglue12pt

At a prime of ordinary reduction, the Iwasawa ``main conjecture''
for elliptic curves relates a Selmer group to a $p$-adic
$L$-function.  In the supersingular case, the statement of the main
conjecture is more complicated as neither the Selmer group nor the
$p$-adic $L$-function is well-behaved.  Recently Kobayashi discovered an
equivalent formulation of the main conjecture at supersingular primes
that is similar in structure to the ordinary case.  Namely, Kobayashi's
conjecture relates modified Selmer groups, which he defined, with modified
$p$-adic $L$-functions defined by the first author.  In this paper we
prove Kobayashi's conjecture for elliptic curves with complex
multiplication.

\vfil
\intro

Iwasawa theory was introduced into the study of the arithmetic of 
elliptic curves by Mazur in the 1970's.  Given an elliptic curve $E$ 
over $\Q$ 
and a prime~$p$ there are two parts to such a program:
an Iwasawa-Selmer module containing information about the arithmetic 
of $E$ over subfields of the cyclotomic\break $\Zp$-extension $\Qinf$ of $\Q$,
and
a $p$-adic $L$-function attached to $E$, belonging to a suitable Iwasawa 
algebra.
The goal, or ``main conjecture'', is to relate these two objects 
by proving that the $p$-adic $L$-function controls 
(in precise terms, is a characteristic power series of the 
Pontrjagin dual of) 
the Iwasawa-Selmer module.  The main conjecture has important 
consequences for the Birch and Swinnerton-Dyer conjecture \pagebreak for $E$.

For primes $p$ where $E$ has ordinary reduction, 
\begin{itemize}
\item
Mazur introduced and studied the Iwasawa-Selmer module \cite{Ma}, 

\item
Mazur and Swinnerton-Dyer constructed the $p$-adic $L$-function \cite{MSD}, 

\item
the main conjecture was proved 
by the second author 
for elliptic curves with complex multiplication \cite{Ru3}, 

\item
Kato proved that the characteristic power series of 
the Pontrjagin dual of the Iwasawa-Selmer module divides the 
$p$-adic $L$-function \cite{Ka}.  
\end{itemize}
The latter two results are proved using Kolyvagin's Euler system machinery.

For primes $p$ where $E$ has supersingular reduction, progress has 
been much slower.  Using the same definitions as for the ordinary case 
gives a Selmer module that is not a torsion Iwasawa module \cite{Ru1}, and a 
$p$-adic $L$-function that does not belong to the Iwasawa algebra \cite{MTT}, \cite{AV}.
Perrin-Riou and Kato made important progress in understanding the case of 
supersingular primes, and independently proposed a main conjecture 
\cite{PR3}, \cite{Ka}.

More recently, the first author \cite{Po} 
proved that when $p$ is a prime of 
supersingular reduction (and either $p > 3$ or $a_p = 0$) 
the ``classical'' $p$-adic $L$-function 
corresponds in a precise way to {\it two} elements $\L_E^+, \L_E^-$ 
of the Iwasawa algebra.  
Shortly thereafter Kobayashi \cite{Ko} defined two submodules 
$\Sel^+(E/\Qinf),\break \Sel^-(E/\Qinf)$ of the ``classical'' Selmer module, 
and proposed a main conjecture: that $\L_E^\pm$ is a characteristic 
power series of the Pontrjagin dual of $\Sel^\pm(E/\Qinf)$.  
Kobayashi proved that this conjecture is equivalent to the 
Kato and Perrin-Riou conjecture, and 
(as an application of Kato's results \cite{Ka}) 
that the characteristic power series of the 
Pontrjagin dual of $\Sel^\pm(E/\Qinf)$ divides~$\L_E^\pm$.

The purpose of the present paper is to prove Kobayashi's main conjecture 
when the elliptic curve $E$ has complex multiplication:  

\vglue8pt {\elevensc Theorem}.
{\it If $E$ is an elliptic curve over $\Q$ with complex multiplication{\rm ,} and $p > 2$ is 
a prime where $E$ has good supersingular reduction{\rm ,} then 
$\L_E^\pm$ is a characteristic power series of the Iwasawa module }
$\Hom(\Sel^\pm(E/\Qinf),\Qp/\Zp)$. 
\vglue8pt

See Definition \ref{kosel} for the definition of Kobayashi's Selmer groups\break
$\Sel^\pm(E/\Qinf)$, and Section~\ref{mcproof} for the definition of $\L_E^\pm$.  With 
the same proof (and a little extra notation) one can prove an analogous 
result for $\Sel^\pm(E/\Q(\bmu_{p^\infty}))$, the Selmer groups over the 
full $p$-cyclotomic field $\Q(\bmu_{p^\infty})$.

The proof relies on the Euler system of elliptic units, and the results and 
methods of \cite{Ru3} which also went into the proof of the 
main conjecture for CM elliptic curves at ordinary primes.  We sketch the ideas 
briefly here, but we defer the precise definitions, statements, and references 
to the main text below.

Fix an elliptic curve $E$ defined over $\Q$ 
with complex multiplication by 
an imaginary quadratic field $K$, and a prime $p > 2$ where $E$ has good 
reduction (ordinary or supersingular, for the moment).  Let $\p$ be a prime 
of $K$ above $p$, and 
let $\K = K(E[p^\infty])$, the (abelian) extension 
of $K$ generated by all $p$-power torsion points on $E$.
Class field theory gives an exact sequence
\begin{equation}
\label{pre4term}
0 \too \E/\C \too \U/\C \too \X \too \A \too 0
\end{equation}
where $\U$, $\E$, and $\C$ are the inverse limits of the local units, 
global units, and elliptic units, respectively, up the tower of 
abelian extensions $K(E[p^n])$ of $K$, and $\X$ (resp.\ $\A$) is the Galois 
group over $K(E[p^\infty])$ of the maximal unramified outside $\p$ 
(resp.\ everywhere unramified) abelian $p$-extension of $K(E[p^\infty])$.  
Further
\vglue4pt
\noindent\hglue6pt\hangindent=24pt\hangafter=1
(a)
the classical Selmer group $\Selgp(E/\K) = \Hom(\X,E[\p^\infty])$, 

\vglue4pt
\noindent\hglue6pt\hangindent=24pt\hangafter=1
(b) 
the ``Coates-Wiles logarithmic derivatives'' of the elliptic units are 
special values of Hecke $L$-functions attached to $E$, 

\vglue4pt
\noindent\hglue6pt\hangindent=24pt\hangafter=1 
(c) 
the Euler system of elliptic units can be used to show that the (torsion) 
Iwasawa modules $\E/\C$ and $\A$ have the same characteristic ideal.
\vglue4pt

If $E$ has ordinary reduction at $p$, then $\U/\C$ and 
$\X$ are torsion Iwasawa modules.  
It then follows from \eqref{pre4term} and (c) that 
$\U/\C$ and $\X$ have the same characteristic ideal, and from 
(b) that the characteristic ideal of $\U/\C$ is a (``two-variable'')
$p$-adic $L$-function.  Now using (a) and restricting to 
$\Qinf \subset \K$ one can prove the main conjecture in this case.

When $E$ has supersingular reduction at $p$, the Iwasawa modules $\U/\C$ and 
$\X$ are {\it not} torsion  (they have rank one), so the argument above breaks down.  
However, Kobayashi's construction suggests a way to remedy this.  Namely, 
one can define submodules $\V^+, \V^- \subset \U$ such that in the exact sequence 
induced from \eqref{pre4term}

\centerline{${\displaystyle
0 \too \E/\C \too \U/(\C+\V^\pm) \too \X/\image(\V^\pm) \too \A \too 0
}$}
\vglue4pt\noindent we have torsion modules $\U/(\C+\V^\pm)$ and $\X/\image(\V^\pm)$, and 
the Kobayashi Selmer groups satisfy
\begin{itemize}
\item[($\hbox{a}'$)]
$\Sel^\pm(E/\Qinf) =  
    \Hom(\X/\image(\V^\pm),E[p^\infty])^{G_{\Qinf}}$. 
\end{itemize}
Using (b) (to relate $\U/(\C+\V^\pm)$ with $\L_E^\pm$)
and (c) as above this will enable us to prove the main 
conjecture in this case as well.

The layout of the paper is as follows.  The general setting and notation 
are laid out in Section~\ref{setup}.  Sections \ref{classsel} and \ref{koselsec} 
describe the classical and Kobayashi Selmer groups, and Sections \ref{kp} and \ref{euerl} 
relate Kobayashi's construction to local units, elliptic units, and $L$-values.  
Section \ref{ru3} applies the results of \cite{Ru3} to our situation.
The proof of the main theorem (restated as Theorem \ref{mainthm} below) 
is given in Section~\ref{mcproof}, and in Section~\ref{appl} we give some arithmetic applications.

\section{The setup}
\label{setup}

Throughout this paper we fix an elliptic curve $E$ defined over $\Q$, with 
complex multiplication by the ring of integers $\O$ of an imaginary quadratic 
field~$K$.  
(No generality is lost by assuming that $\End(E)$  is the maximal order in~$K$, 
since we could always replace $E$ by an isogenous curve with this property.)
Fix also a rational prime $p > 2$ where $E$ has good supersingular 
reduction.  As is well known, it follows that $p$ remains prime in $K$.  
It also follows that $a_p = p+1-|E(\Fp)| = 0$, so 
we can apply the results of the first author \cite{Po} 
and Kobayashi \cite{Ko}.
Let $\Kp$ and $\Op$ denote the completions of $K$ and $\O$ at $p$.

For every $k$ let $E[p^k]$ denote kernel of $p^k$ in $E(\bar{\Q})$, 
$E[p^\infty] = \mathbold{\cup}_k E[p^k]$, and $T_p(E) = \varprojlim E[p^k]$.  
Let $\K = K(E[p^\infty])$, let $\Kinf$ 
denote the (unique)\break $\Zp^2$-extension of $K$, let $\Qinf \subset \Kinf$ be the 
cyclotomic $\Zp$-extension of $\Q$, and let $\KQinf = K\Qinf \subset \Kinf$ 
be the cyclotomic $\Zp$-extension of $K$.  
Let $\rho$ denote the character
$$
\rho : G_K \too \Aut_{\Op}(E[p^\infty]) \cong \Opx.
$$
Let $\Efg$ denote the formal group giving the kernel of reduction modulo $p$ on~$E$.  
The theory of complex multiplication shows that $\Efg$ 
is a Lubin-Tate formal group of height two over $\Op$ for the uniformizing 
parameter $-p$.  It follows that 
$\rho$ is surjective, even when restricted to an inertia 
group of $p$ in $G_K$.  Therefore  
$p$ is totally ramified in $\K/K$ and $\rho$ induces an isomorphism
$\Gal(\K/K) \cong \Opx$.  We can decompose
$$
\Gal(\K/K) = \Delta \times \Gamma_+ \times \Gamma_-
$$
where $\Delta = \Gal(\K/\Kinf) \cong \Gal(K(E[p])/K)$ is the non-$p$ part of 
$\Gal(\K/K)$, which is 
cyclic of order $p^2-1$, and $\Gamma_\pm$ is the largest subgroup 
of $\Gal(\K/K(E[p]))$ on which the nontrivial element of $\Gal(K/\Q)$ acts by 
$\pm 1$.  Then $\Gamma_+$ and $\Gamma_-$ are both free of rank one over 
$\Zp$.   

Let $M$ (resp.\ $L$) denote the maximal abelian $p$-extension of 
$K(E[p^\infty])$ that is unramified outside of the unique prime above $p$ 
(resp.\ unramified\break everywhere), and let $\X = \Gal(M/\K)$ and 
$\A = \Gal(L/\K)$.  If $F$ is a finite extension of $K$ in $\K$ let 
$\O_{F}$ denote the ring of integers of $F$, and define subgroups 
$C_F \subset E_F \subset U_F \subset (\O_F \otimes \Zp)^\times$ as follows.  
The group $U_F$ is the pro-$p$-part of the local unit group 
$(\O_F \otimes \Zp)^\times$,
$E_F$ is the closure of the projection of the global units $\O_F^\times$  
into $U_F$, and $C_F$ is the closure of the projection of the subgroup of 
elliptic units (as defined for example in \S1 of \cite{Ru3}) into $U_F$.  
Finally, define 
$$
\C = \varprojlim C_F \subset \E = \varprojlim E_F \subset \U = \varprojlim U_F,  
$$
inverse limit \pagebreak with respect to the norm map over finite extensions of $K$ in $\K$.
Class field theory gives an isomorphism $\Gal(M/L) \cong \U/\E$.  
We summarize this setting in Figure 1 below.
\figin{rubin}{1000}
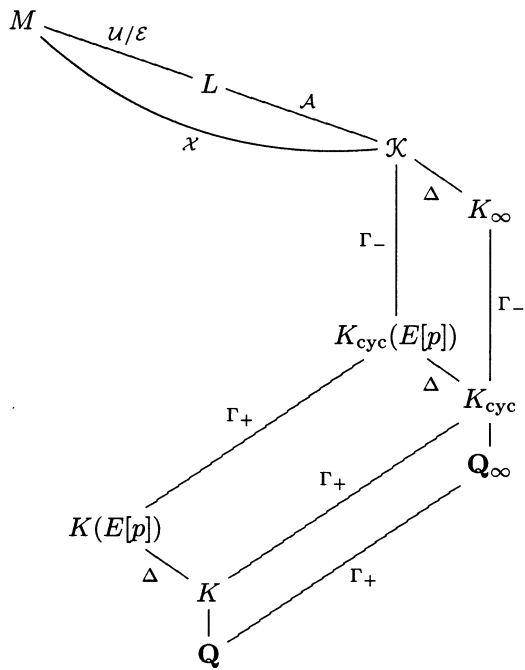
\centerline{Figure 1.}
\vglue12pt

If $K \subset F \subset \K$ we define the Iwasawa algebra $\Iw{F} = \Zp[[\Gal(F/K)]]$.
In particular we have
\begin{eqnarray*}
\Iw{\K} &=& \Zp[[\Gal(\K/K)]] = \Zp[[\Delta \times \Gamma_+ \times \Gamma_-]], \\ 
\Iw{\Kinf}& = &\Zp[[\Gal(\Kinf/K)]] = \Zp[[\Gamma_+ \times \Gamma_-]], \\ 
\Iw{\KQinf} &=& \Zp[[\Gal(\KQinf/K)]] \cong \Zp[[\Gamma_+]] \cong \Zp[[\Gal(\Qinf/\Q)]].
\end{eqnarray*}
We write simply $\Iwq$ for $\Iw{\KQinf}$, and we write $\IwO{F} = \Iw{F} \otimes \Op$ 
and $\IwOq = \Iwq \otimes \Op$.

\numbereddemo{Definition}
Suppose $Y$ is a $\Iw{\K}$-module.  We define the twist
$$
\tw{Y} = Y \otimes \Hom_\O(E[p^\infty],\Kp/\Op).
$$ 
The module $\Hom_\O(E[p^\infty],\Kp/\Op)$ is free of rank one over $\Op$, 
and $G_K$ acts on it via $\rho^{-1}$.  Thus we have 
$\tw{T_p(E)} \cong \Op$ and $\tw{E[p^\infty]} \cong \Kp/\Op$.

If $K \subset F \subset \K$ we define
$$
\redgen{F}{Y} = \tw{Y} \otimes_{\Iw{\ssK}} \Iw{F} 
    = \tw{Y}/\langle \gamma-1 : \gamma \in \Gal(\K/F)\rangle,
$$
the $F$-coinvariants of $\tw{Y}$.
We will be interested in $\rc{Y}$ and $\red{Y}$.
Concretely, if we write $Z$ for the $\IwO{\Kinf}$-submodule of 
$Y \otimes \Op$ on which $\Delta$ acts via $\rho$, then 
$\rc{Y}$ can be identified with $\tw{Z}$ and 
$\red{Y}$ can be identified with $\tw{(Z/(\gaminus - \rho(\gaminus))Z)}$ where 
$\gaminus$ is a topological generator of $\Gamma_-$.
\enddemo

\section{The classical Selmer group}
\label{classsel}

For every number field $F$ we have the classical $p$-power Selmer group 
$\Sel(E/F) \subset H^1(F,E[p^\infty])$, which sits in an exact sequence 
$$
0 \too E(F) \otimes (\Qp/\Zp) \too \Sel(E/F) \too \Sh(E/F)[p^\infty] \too 0
$$
where $\Sh(E/F)[p^\infty]$ is the $p$-part of the Tate-Shafarevich group 
of $E$ over $F$.  Taking direct limits allows us to define $\Sel(E/F)$ 
for every algebraic extension $F$ of $\Q$.

\proclaim{Theorem}
\label{selthm}
$\Sel(E/\KQinf) \cong \Hom_\O(\red{\X},\Kp/\Op)$.
\endproclaim

\demo{Proof}
Combining Theorem 2.1, Proposition 1.1, and Proposition 1.2 of \cite{Ru1} shows 
that
\begin{eqnarray*}
\Sel(E/\KQinf) &\cong& \Hom_\O(\X,E[p^\infty])^{\Gal(\ssK/\KQinf)} \\
    &=& \Hom_\O(\tw{\X},\Kp/\Op)^{\Gal(\ssK/\KQinf)} = \Hom_\O(\red{\X},\Kp/\Op).\\
\noalign{\vskip-24pt}
\end{eqnarray*}
\enddemo
\vglue9pt

\numbereddemo{{R}emark}
We have $\rank_{\IwO{\Kinf}} \rc{\X} = 1$ (see for example [Ru3,\break Th.\ 5.3(iii)]), so
$\rank_{\IwOq}\red{\X} \ge 1$. Thus, unlike the case of ordinary primes, the Selmer group 
$\Sel(E/\KQinf)$ is not a co-torsion $\IwOq$-module.  
This makes the Iwasawa theory for supersingular primes more difficult 
than the ordinary case.  In the next section, following Kobayashi \cite{Ko}, 
we will remedy this by defining two smaller Selmer groups which will both 
be co-torsion $\IwOq$-modules.
\enddemo

\section{Kobayashi's restricted Selmer groups}
\label{koselsec}
\advance\eqcount by 1

If $F$ is a finite extension of $K$ in $\K$ let $F_p$ denote the completion of $F$ 
at the unique prime above $p$, and for an arbitrary $F$ with $K \subset F \subset \K$ 
let $F_p = \mathbold{\cup}_{N}N_p$, union over finite extensions of $K$ in $F$.  For every such $F$
let $\m_F$ denote the maximal ideal of $F_p$ and   
let $E_1(F_p) \subset E(F_p)$ be the kernel of reduction.  Then $E_1(F_p)$ 
is the pro-$p$ part of $E(F_p)$ and we define the logarithm map $\logE$ to be the 
composition 
$$
\logE : E(F_p) \onto E_1(F_p) \isom \Efg(\m_F) \too F_p
$$
where the first map is projection onto the pro-$p$ part, the second is the canonical 
isomorphism between the kernel of reduction and the formal group $\Efg$, and the third 
is the formal group logarithm map.

\numbereddemo{Definition}
For $n \ge 0$ let $\Qn$ denote the extension of $\Q$ of degree $p^n$ in $\Qinf$, 
and if $n \ge m$ let $\Tr_{n/m}$ 
denote the trace map from $E(\Qnp)$ to $E(\Q_{m,p})$.  
For each $n$ define two subgroups $E^+(\Qnp), E^-(\Qnp) \subset E(\Qnp)$ by
\begin{eqnarray*}
E^+(\Qnp) &= &\{x \in E(\Qnp) : \hbox{$\Tr_{n/m}x \in E(\Q_{m-1,p})$ if $0 < m \le n$, $m$ odd}\} \\
E^-(\Qnp) &= &\{x \in E(\Qnp) : \hbox{$\Tr_{n/m}x \in E(\Q_{m-1,p})$ if $0 < m \le n$, $m$ even}\}
\end{eqnarray*}
and let $E_1^\pm(\Qnp) = E^\pm(\Qnp) \cap E_1(\Qnp)$. 
Equivalently, let $\Xi_n^+$ (resp.\ $\Xi_n^-$) denote the 
set of nontrivial characters  
$\Gal(\Qn/\Q) \to \bmu_{p^n}$ whose order is an odd (resp.\ even) 
power of $p$, and then 
$$
E^\pm(\Qnp) = \{x \in E(\Qnp) : 
	\hbox{$\sum_{\sigma \in \Gal(\Qn/\Q)} \chi(\sigma) x^\sigma = 0$ for every $\chi \in \Xi_n^\pm$}\}
$$
where the sum takes place in $E(\Qnp) \otimes \Z[\bmu_{p^n}]$.  
Note that when $n = 1$ we get $E^+(\Qp) = E^-(\Qp) = E(\Qp)$.
When $n = \infty$ we define
$$
E^\pm(\Qinfp) = \mathbold{\cup}_n E^\pm(\Qnp).
$$
We also define $E^\pm(K\Qnp)$ exactly as above with $\Qn$ replaced by $K\Qn$.  The 
complex multiplication map $E(\Qnp) \otimes \Op \to E(K\Qnp)$ induces isomorphisms 
\begin{equation}
\label{qtok}
E_1(\Qnp) \otimes \Op \isom E_1(K\Qnp), 
    \quad E^\pm_1(\Qnp) \otimes \Op \isom E^\pm_1(K\Qnp)
\end{equation}
for every $n \le \infty$.
\enddemo

Fix once and for all a generator $\{\zeta_{p^n}\}$ of $\Zp(1)$, so $\zeta_{p^n}$ is 
a primitive $p^n$-th root of unity and $\zeta_{p^{n+1}}^p = \zeta_{p^n}$.  If 
$\chi : \Gamma_+ \onto \bmu_{p^k}$ define the Gauss sum 
$\tau(\chi) = \sum_{\sigma \in \Gal(\Q(\bmu_{p^k})/\Q)} \chi(\sigma)\zeta_{p^k}^\sigma$.

\proclaimtitle{Kobayashi \cite{Ko}}
\proclaim{Theorem}
\label{kob}
\begin{itemize}
\ritem{(i)}
$E^+(\Qnp) + E^-(\Qnp) = E(\Qnp)$.
\ritem{(ii)}
$E^+(\Qnp) \cap E^-(\Qnp) = E(\Qp)$.
\end{itemize}
Further{\rm ,} there is a sequence of points $d_n \in E_1(\Qnp)$ {\rm (}\/depending on the choice 
of $\{\zeta_{p^n}\}$ above\/{\rm )} with the following properties.
\begin{itemize}
\ritem{(iii)}
$\Tr_{n/n-1} d_n = 
\left\{ \begin{array}{ll}
d_{n-2} & \hbox{if $n\ge 2$}, \\
\frac{1-p}{2}d_0 & \hbox{if $n = 1$.}
\end{array}\right. $
\ritem{(iv)}
If $\chi : \Gal(\Qn/\Q) \isom \bmu_{p^n}$ then
$$
\sum_{\sigma \in \Gal(\Qn/\Q)} \chi(\sigma) \logE(d_n^\sigma) = 
\left\{ \begin{array}{ll}
(-1)^{[\frac{n}{2}]}\tau(\chi) & \hbox{if $n > 0$}, \\
\frac{p}{p+1} & \hbox{if $n = 0$.}
\end{array}\right. 
$$
\ritem{(v)}
If $\varepsilon = (-1)^n$ then 
\end{itemize}
\hfill${\displaystyle
E_1^\varepsilon(\Qnp) = \Zp[\Gal(\Qn/\Q)]d_n \quad\hbox{and}\quad
E_1^{-\varepsilon}(\Qnp) = \Zp[\Gal(\Q_{n-1}/\Q)]d_{n-1}.
}$\hfill
\endproclaim
 
\vglue6pt
{\it Proof}.
The first two assertions are Proposition 8.12(ii) of \cite{Ko}.

Let $d_n = (-1)^{[\frac{n+1}{2}]}\Tr_{\Q(\bmu_{p^{n+1}})/\Qn}c'_{n+1}$ where 
$c'_{n+1} \in E_1(\Q(\bmu_{p^{n+1}})_p)$ corresponds to the point 
$c_{n+1} \in \Efg(\Q(\bmu_{p^{n+1}})_p)$ defined by Kobayashi in Section~4 of \cite{Ko}.
Then the last three assertions of the theorem follow from 
Lemma 8.9, Proposition 8.26, and Proposition 8.12(i), 
respectively, of \cite{Ko}.
\hfill\qed

\numbereddemo{Definition}
\label{kosel}
If $0 \le n \le \infty$ we define Kobayashi's restricted Selmer groups 
$\Sel^\pm(E/\Qn) \subset \Sel(E/\Qn)$ by
$$
\Sel^\pm(E/\Qn) = \ker \Bigl(\Sel(E/\Qn) \to 
    H^1(\Qnp,E[p^\infty])/(E^\pm(\Qnp) \otimes \Qp/\Zp)\Bigr).
$$
Since $E(\Q_{n,v}) \otimes \Qp/\Zp = 0$ when $v \nmid p$, 
a class $c \in H^1(\Qn,E[p^\infty])$ belongs to $\Sel^\pm(E/\Qn)$ 
if and only if its localizations $c_v \in H^1(\Q_{n,v},E[p^\infty])$ satisfy 
$c_v = 0$ if $v \nmid p$ and
$$
c_p \in \image \Bigl(E^\pm(\Qnp) \otimes \Qp/\Zp \to H^1(\Qnp,E[p^\infty])\Bigr).
$$
(If we replace $E^\pm(\Qnp)$ by $E(\Qnp)$ we get the definition of $\Sel(E/\Qn)$.)

We define $\Sel^\pm(E/\KQinf)$ in exactly the same way with $\Qn$ replaced by $K\Qn$, 
using $E^\pm(K\Qnp)$, and then
$$
\Sel^\pm(E/\Qinf) \otimes \Op \cong \Sel^\pm(E/\KQinf).
$$
\enddemo

\vglue-16pt
\section{The Kummer pairing}
\label{kp}
\advance\eqcount by 2

The composition
\begin{eqnarray*}
&&\hskip-.5in E(\K_p) \otimes \Qp/\Zp \too H^1(\K_p,E[p^\infty]) \isom \Hom(G_{\ssK_p},E[p^\infty]) \\
  &&\hskip.75in  \too \Hom(\U,E[p^\infty]) \isom \Hom_\O(\tw{\U},\Kp/\Op),
\end{eqnarray*}
where the third map is induced by the inclusion $\U \hookrightarrow G_{\ssK_p}$ 
of local class field theory, induces an $\Op$-linear Kummer pairing 
\begin{equation}
\label{pairing}
(E(\K_p) \otimes \Qp/\Zp) \times \tw{\U} \to \Kp/\Op.
\end{equation}

\proclaim{Proposition}
\label{lulpprop}
The Kummer pairing of {\rm \eqref{pairing}} induces an isomorphism
\vglue12pt
\hfill ${\displaystyle
\red{\U} \cong \Hom_\O(E(\KQinfp) \otimes \Qp/\Zp,\Kp/\Op). 
}$\hfill
\endproclaim

{\it Proof}.
This is equivalent to Proposition 5.4 of \cite{Ru2}.
\hfill\qed

\numbereddemo{Definition}
Define $\redd{\V}^\pm \subset \red{\U}$ to be the subgroup of $\red{\U}$ 
corresponding to $\Hom_\O(E(\KQinfp)/E^\pm(\KQinfp) \otimes \Qp/\Zp,\Kp/\Op)$ 
under the isomorphism of Proposition \ref{lulpprop}.  Since $\Hom_\O(~\cdot~,\Kp/\Op)$ 
is an exact functor on $\Op$-modules we have 
\begin{eqnarray}
\label{lplu}
&&E^\pm(\KQinfp) \otimes \Qp/\Zp \cong \Hom_\O(\red{\U}/\redd{\V}^\pm,\Kp/\Op), \\
\label{lulp}
&&\red{\U}/\redd{\V}^\pm \cong \Hom_\O(E^\pm(\KQinfp) \otimes \Qp/\Zp,\Kp/\Op).
\end{eqnarray}
\enddemo

Let $\alpha : \U \to \X$ be the Artin map of global class field theory.  
The following theorem is the step labeled ($\hbox{a}'$) in the introduction.

\proclaim{Theorem}
\label{newsel}
$\Sel^\pm(E/\KQinf) = \Hom_\O(\red{\X}/\alpha(\redd{\V}^\pm),\Kp/\Op)$.
\endproclaim

{\it Proof}.
This is Theorem \ref{selthm} combined with Definition \ref{kosel} of $\Sel^\pm(E/\KQinf)$ 
and \eqref{lplu}.
\hfill\qed

\proclaim{Proposition}
\label{ufree}
{\rm (i)}
$\rc{\U}$ is free of rank two over $\IwO{\Kinf}$ and $\red{\U}$ is free of rank two 
over $\IwOq$.
\begin{itemize}
\ritem{(ii)}
$\redd{\V}^\pm$ and $\red{\U}/\redd{\V}^\pm$ are free of rank one over $\IwOq$.
\ritem{(iii)}
There is a {\rm (}\/noncanonical\/{\rm )} submodule $\V^\pm \subset \rc{\U}$ whose image in  
$\red{\U}$ is $\redd{\V}^\pm$ and such that $\V^\pm$ and $\rc{\U}/\V^\pm$ 
are free of rank one over $\IwO{\Kinf}$.
\end{itemize}

\endproclaim

{\it Proof}.
By \cite{Gr}, $\rc{\U}$ is free of rank two over $\IwO{\Kinf}$, and then 
the definition of $\red{\U}$ shows that $\red{\U}$ is free of rank two over $\IwOq$.
Theorem 6.2 of \cite{Ko} (see also Theorem \ref{genmu} below)
and \eqref{lulp} show that $\red{\U}/\redd{\V}^\pm$ 
is free of rank one over $\IwOq$, so the exact sequence 
$0 \to \redd{\V}^\pm \to \red{\U} \to \red{\U}/\redd{\V}^\pm \to 0$ 
splits.  Thus $\redd{\V}^\pm$ is a projective $\IwOq$-module, and 
Nakayama's lemma shows that every projective $\IwOq$-module is free.  
This proves (ii).

Let $u$ be any element of $\rc{\U}$ whose image in $\red{\U}$ generates 
$\redd{\V}^\pm$, and let $\V^\pm = \IwO{\Kinf}u$.  Then $\V^\pm$ is free of 
rank one, and it follows from (ii) and Nakayama's lemma that 
$\rc{\U}/\V^\pm$ is free of rank one over $\IwO{\Kinf}$ as well. \phantom{move}
\hfill\qed
 
\vglue-6pt
\section{Elliptic units and the explicit reciprocity law}
\label{euerl}
\vglue-6pt

Let $\psi_E$ denote the Hecke character of $K$ attached to $E$, and for 
every character $\chi$ of finite order of $G_K$ let $L(\psi_E\chi,s)$ denote 
the Hecke $L$-function.  If $\chi$ is the restriction of a character of $G_\Q$ 
then $L(\psi_E\chi,s) = L(E,\chi,s)$, the usual $L$-function of $E$ twisted 
by the Dirichlet character $\chi$.  Let $\Omega_E \in \R^+$ denote the 
real period of a minimal model of $E$. 

The explicit reciprocity law of Wiles \cite{Wi} together with 
a computation of Coates and Wiles \cite{CW} leads to the 
following theorem, which is the step labeled (b) in the introduction.

\proclaim{Theorem}
\label{erl}
The $\IwO{\Kinf}$\/{\rm -}\/module $\rc{\C}$ of elliptic units is free of rank one 
over $\IwO{\Kinf}$.  It has a generator $\xi$ with the property that 
if $K \subset F \subset \Kinf${\rm ,} 
$x \in E(F_p)${\rm ,} and $\chi : \Gal(F/K) \to \bmu_{p^\infty}${\rm ,} 
then the Kummer pairing $\ld ~,~ \rd$ of {\rm \eqref{pairing}} satisfies
$$
\sum_{\sigma\in\Gal(F/K)} \chi^{-1}(\sigma) \ld x^\sigma \otimes p^{-k}, \xi \rd 
    = p^{-k} \frac{L(\psi_E\chi,1)}{\Omega_E} 
        \sum_{\sigma\in\Gal(F/K)}\chi^{-1}(\sigma)\lambda_E(x^\sigma).
$$
\endproclaim

{\it Proof}.
See \cite{Wi} and \cite[\S5]{CW}, or Theorem 7.7(i) of 
\cite{Ru3} and Theorem~3.2 and the proof of Proposition 5.6 of \cite{Ru2}.
\hfill\qed

\proclaim{{C}orollary}
\label{cnotinu}
{\rm (i)} The map $\red{\C} \to \red{\U}$ is injective.
\begin{itemize}
\ritem{(ii)}
$\red{\C}$ is free of rank one over $\IwOq$ and 
$\red{\C} \cap \redd{\V}^+ = \red{\C} \cap \redd{\V}^- =  0$.
\ritem{(iii)}
$\rank_{\IwO{\Kinf}}\rc{\E} = 1$ and
$\rc{\E} \cap \V^+ = \rc{\E} \cap \V^- =  0$.
\end{itemize}

\endproclaim

{\it Proof}.
Since $\red{\C}$ and $\red{\U}/\redd{\V}^\pm$ are free of rank one 
over $\IwOq$ (Theorem \ref{erl} and Proposition \ref{ufree}(ii)), 
the map $\red{\C} \to \red{\U}/\redd{\V}^\pm$ is 
either injective or identically 
zero.  Thus to prove both (i) and (ii) it will suffice to show that 
the image $\redd{\xi} \in \red{\U}$ of the generator $\xi \in \rc{\C}$ 
of Theorem \ref{erl} satisfies
$\redd{\xi} \notin \redd{\V}^+$ and $\redd{\xi} \notin \redd{\V}^-$.

Rohrlich \cite{Ro} proved that $L(E,\chi,1) \ne 0$ 
for all but finitely many characters $\chi$ of $\Gal(\KQinf/K)$.  
Applying Theorem \ref{erl} with $x = d_{2n}$ for large $n$ 
and using Theorem \ref{kob}(iv) it follows that 
the image of $\xi$ in $\Hom_\O(E^+(\KQinfp) \otimes \Qp/\Zp,\Kp/\Op)$ 
is nonzero.  Hence $\redd{\xi} \notin \redd{\V}^+$.  Similarly, using the points 
$d_{2n+1}$ for large $n$ shows that $\redd{\xi} \notin \redd{\V}^-$.  
This proves (i) and (ii).

By Corollary 7.8 of \cite{Ru3}, $\rc{\E}$ is a torsion-free, rank-one 
$\IwO{\Kinf}$-module.  Just as in (i), since $\rc{\U}/\V^\pm$ is torsion-free 
(Proposition \ref{ufree}(iii)) the map $\rc{\E} \to \rc{\U}/\V^\pm$ is either 
injective or identically zero.  But we saw above that 
$\redd{\xi} \notin \redd{\V}^\pm$, so $\xi \notin \V^\pm$ 
and $\rc{\E} \to \rc{\U}/\V^\pm$ is not identically zero.  This proves (iii).
\hfill\qed

\section{The characteristic ideals}
\label{ru3}
\advance\eqcount by 5

If $B$ is a finitely generated torsion module over $\IwO{\Kinf}$ 
(resp.\ $\IwOq$,\break resp.\ $\Iwq$), we will write $\cp_{\IwO{\Kinf}}(B)$ 
(resp.\ $\cp_{\IwOq}(B)$, resp.\ $\cp_{\Iwq}(B)$) for its characteristic 
ideal.

The following theorem is Theorem 4.1(ii) of \cite{Ru3}, twisted by $\rho^{-1}$.  
It is the step labeled (c) in the introduction.

\proclaimtitle{\cite{Ru3}}
\proclaim{Theorem}
\label{oldmc}
The $\IwO{\Kinf}$\/{\rm -}\/modules $\rc{\A}$ and $\rc{\E}/\rc{\C}$ 
are fi\-nitely generated and torsion{\rm ,} and
\vglue9pt \hfill${\displaystyle
\cp_{\IwO{\Kinf}}(\rc{\A}) = \cp_{\IwO{\Kinf}}(\rc{\E}/\rc{\C}).
}$\hfill
\endproclaim
\vglue4pt
 
{\elevensc {C}orollary 6.2.}
{\it Let $\alpha : \U \to \X$ denote the Artin map of global class field theory.  Then 
$\rc{\X}/\alpha(\V^\pm)$ and $\rc{\U}/(\V^\pm+\rc{\C})$ are finitely generated 
torsion $\IwO{\Kinf}$\/{\rm -}\/modules and}
$$
\cp_{\IwO{\Kinf}}(\rc{\X}/\alpha(\V^\pm)) 
    = \cp_{\IwO{\Kinf}}(\rc{\U}/(\V^\pm+\rc{\C})).
$$
\vglue3pt
\advance\theoremcount by 1

{\it Proof}.
Class field theory gives an exact sequence
$$
0 \too \E/\C \too \U/\C \map{~\alpha~} \X \too \A \to 0.
$$
Twisting by $\rho^{-1}$ and using the fact that $\Delta$ has order 
prime to $p$ gives another exact sequence
$$
0 \too \rc{\E}/\rc{\C} \too \rc{\U}/\rc{\C} \map{~\alpha~} 
    \rc{\X} \too \rc{\A} \too 0.
$$
Since $\rc{\E} \cap \V^\pm = 0$ by Corollary \ref{cnotinu}, we get 
finally an exact sequence
\begin{equation}
\label{modseq}
0 \to \rc{\E}/\rc{\C} \to \rc{\U}/(\V^\pm+\rc{\C}) \map{\alpha} 
    \rc{\X}/\alpha(\V^\pm) \to \rc{\A} \to 0.
\end{equation}
Since $\rc{\C} \cap \V^\pm = 0$,  
it follows from Theorem \ref{erl} and Proposition \ref{ufree} that 
the quotient 
$\rc{\U}/(\V^\pm+\rc{\C})$ is a finitely generated torsion 
$\IwO{\Kinf}$-module.  Now \eqref{modseq} and Theorem \ref{oldmc}
show that $\rc{\X}/\alpha(\V^\pm)$ is a finitely generated torsion 
$\IwO{\Kinf}$-module as well, and that the two characteristic ideals 
are equal.
\hfill\qed

\proclaim{Theorem}
\label{bettermc}
The $\IwOq$\/{\rm -}\/modules $\red{\X}/\alpha(\redd{\V}^\pm)$ and 
$\red{\U}/(\redd{\V}^\pm + \red{\C})$ are fi\-nite\-ly generated 
torsion modules and
$$
\cp_{\IwOq}(\red{\X}/\alpha(\redd{\V}^\pm)) 
    = \cp_{\IwOq}(\red{\U}/(\redd{\V}^\pm + \red{\C})).
$$
Further{\rm ,} $\red{\X}/\alpha(\redd{\V}^\pm)$ has no finite $\IwOq$\/{\rm -}\/submodules.
\endproclaim

The proof of Theorem \ref{bettermc} is given below, 
after a few lemmas.  The proof is essentially contained in Section~11 of \cite{Ru3}, 
but since it is crucial for our main result we recall some of the details.

If $\AA$ is an ideal of $\IwO{\Kinf}$, let $\overline{\AA} \subset \IwOq$ 
denote the image of $\AA$ under the projection map $\IwO{\Kinf} \onto \IwOq$.  
Fix a topological generator $\gaminus$ of $\Gamma_- = \Gal(\Kinf/\KQinf)$.

\proclaim{Lemma}
\label{redcp}
Suppose $B$ is a finitely generated torsion $\IwO{\Kinf}$\/{\rm -}\/module with no nonzero 
pseudo\/{\rm -}\/null submodules.
Then 
$$
\hbox{$\overline{\cp_{\IwO{\Kinf}}(B)} \ne 0$ if and only if $B/(\gaminus-1)B$ 
is a torsion $\IwOq$\/{\rm -}\/module}\/, 
$$
and in that case 
$$
\cp_{\IwOq}(B/(\gaminus-1)B) = \overline{\cp_{\IwO{\Kinf}}(B)}.
$$
\endproclaim

{\it Proof}.
See Lemma 4 of  \cite[\S{I.1.3}]{PR1} or Lemma 6.2 of \cite{Ru3}.
\hfill\qed

\proclaim{Lemma}
\label{nopnsub}
Suppose $B$ is a finitely generated $\IwO{\Kinf}$\/{\rm -}\/module with no nonzero 
pseudo\/{\rm -}\/null submodules.  If $B'$ is a free $\IwO{\Kinf}$\/{\rm -}\/submodule of $B$ 
then $B/B'$ has no nonzero pseudo\/{\rm -}\/null submodules.
\endproclaim

{\it Proof}.
By induction we may reduce to the case that $B'$ is free of rank one, 
and may reduce further to the case that $B/B'$ is pseudo-null.  Since 
$\IwO{\Kinf}$ is a unique factorization domain it follows that $B = B'$.
\hfill\qed

\proclaim{Lemma}
\label{morenopnsub}
Suppose $B$ is a finitely generated torsion $\IwO{\Kinf}$\/{\rm -}\/module with no 
nonzero pseudo\/{\rm -}\/null submodules{\rm ,} and both $B/(\gaminus-1)B$ and\break
$B/(\gaminus-\rho^{-1}(\gaminus))B$ are torsion $\IwOq$\/{\rm -}\/modules.  Then 
$B/(\gaminus-1)B$ has a nonzero finite submodule
if and only if $B/(\gaminus-\rho^{-1}(\gaminus))B$ has.
\endproclaim

{\it Proof}.
This is Lemma 11.15 of \cite{Ru3}
\hfill\qed

\demo{Proof of Theorem {\rm \ref{bettermc}}}
By Proposition \ref{ufree} and Corollary \ref{cnotinu}, 
$\red{\U}$ and $\redd{\V}^\pm + \red{\C}$ are free of rank two over $\IwOq$, and  
$\rc{\U}$ and $\V^\pm + \rc{\C}$ are free of rank two over $\IwO{\Kinf}$.
Therefore (using Lemma \ref{nopnsub}) 
$\red{\U}/(\redd{\V}^\pm + \red{\C})$ and $\rc{\U}/(\V^\pm + \rc{\C})$
are torsion modules with no nonzero pseudo-null submodules.  By Lemma \ref{redcp} 
it follows that 
\begin{equation}
\label{cpu}
\cp_{\IwOq}(\red{\U}/(\redd{\V}^\pm + \red{\C})) 
    = \overline{\cp_{\IwO{\Kinf}}(\rc{\U}/(\V^\pm+\rc{\C}))} \ne 0.
\end{equation}

Class field theory shows that the kernel of $\alpha: \rc{\U} \to \rc{\X}$ is $\rc{\E}$.  
Therefore by Corollary \ref{cnotinu} $\alpha$ is injective on $\V^\pm$, so 
$\alpha(\V^\pm)$ is a free, rank-one $\IwO{\Kinf}$-submodule of $\rc{\X}$. 
By \cite{Gr}, $\rank_{\IwO{\Kinf}}\rc{\X} = 1$ and $\rc{\X}$ has no nonzero pseudo-null 
submodules, so (using Lemma \ref{nopnsub}) $\rc{\X}/\alpha(\V^\pm)$ is a torsion 
$\IwO{\Kinf}$-module with no nonzero pseudo-null submodules.  
Further, Corollary~6.2 and \eqref{cpu} show that 
\begin{equation}
\label{cpx}
\overline{\cp_{\IwO{\Kinf}}(\rc{\X}/\alpha(\V^\pm))} 
   = \overline{\cp_{\IwO{\Kinf}}(\rc{\U}/(\V^\pm+\rc{\C}))} \ne 0.
\end{equation}
Thus we can apply Lemma \ref{redcp} to conclude that 
$$
\cp_{\IwOq}(\red{\X}/\alpha(\redd{\V}^\pm)) 
    = \overline{\cp_{\IwO{\Kinf}}(\rc{\X}/\alpha(\V^\pm))},
$$
and together with \eqref{cpu} and \eqref{cpx} this proves 
$$
\cp_{\IwOq}(\red{\X}/\alpha(\redd{\V}^\pm)) 
    = \cp_{\IwOq}(\red{\U}/(\redd{\V}^\pm + \red{\C})).
$$

It remains to prove that $\red{\X}/\alpha(\redd{\V}^\pm)$ has no nonzero 
finite submodules.  This will follow from Lemma \ref{morenopnsub}.  
We give the argument briefly here; 
see the proof of Theorem 11.16 of \cite{Ru3} for more details.

We can identify $\rc{\X}/(\gaminus-\rho^{-1}(\gaminus))\rc{\X}$ 
with a subgroup of $$\tw{(\X/(\gaminus-1)\X)}.$$ 
Standard techniques (for example \cite[\S2]{Gr}) identify $\X/(\gaminus-1)\X$ 
with a subgroup of $\Gal(M_0/\KQinf(E[p]))$ where $M_0$ is the 
maximal abelian $p$-extension of $\KQinf(E[p])$ unramified outside $p$, and 
by \cite{Gr}, $\Gal(M_0/\KQinf(E[p]))$ has no nonzero finite 
submodules.  Hence $\rc{\X}/(\gaminus-\rho^{-1}(\gaminus))\rc{\X}$ 
has no nonzero finite submodules.

Let $B = \rc{\X}/\alpha(\V_\pm)$.  Lemma \ref{nopnsub} now shows that 
$B/(\gaminus-\rho^{-1}(\gaminus))B$ has no nonzero finite submodules, 
and we observed above that $B$ has no nonzero pseudo-null submodules, so 
Lemma \ref{morenopnsub} shows that 
$B/(\gaminus-1)B = \red{\X}/\alpha(\redd{\V}^\pm)$ 
has no nonzero finite submodules.  
\enddemo

 \vglue-9pt
\section{Local units, elliptic units, and the $p$-adic $L$-functions}
\label{mcproof}
\advance\eqcount by 8

Fix a topological generator $\gamma$ of 
$\Gamma_+ \cong \Gal(\KQinf/K) \cong \Gal(\Qinf/\Q)$.  For every $n \ge 1$ define
$$
\nu_n = \sum_{i=0}^{p-1} \gamma^{ip^{n-1}} \in \Iwq
$$
and define $\omega_n^\pm \in \Iwq$ by
$$
\omega^+_n = \prod_{1 \le i \le n,2 \mid i} \nu_i, \quad
    \omega^-_n = \prod_{1 \le i \le n,2 \nmid i} \nu_i.
$$

\proclaimtitle{Kobayashi \cite{Ko}}
\proclaim{Theorem}
\label{genmu}
The $\IwOq$\/{\rm -}\/module 
$$\Hom(E^\pm(\Qinfp) \otimes \Qp/\Zp, \Qp/\Zp)$$ is free of rank one{\rm ,} 
with a generator $\mu^\pm$ such that for every $k, n \in \Z^+${\rm ,} and 
every character $\chi : \Gal(\Qn/\Q) \to \bmu_{p^n}${\rm ,} 
$$
\sum_{\sigma\in\Gal(\Qn/\Q)} \chi(\sigma) \mu^\pm(d_n^\sigma \otimes p^{-k}) 
    = \chi(\omega_n^\mp)p^{-k}.
$$
\endproclaim

{\it Proof}.
An easy exercise shows that for $0 \le n \le \infty$
\begin{equation}
\label{compact}
\Hom(E^\pm(\Qnp) \otimes \Qp/\Zp, \Qp/\Zp) = \Hom(E^\pm(\Qnp),\Zp).
\end{equation}
In Section~8 of \cite{Ko}, especially Proposition 8.18 and Theorem 6.2, 
Kobayashi shows that for every $n$ and $\varepsilon = \pm 1$, the map
$$
f \mapsto 
\left\{ \begin{array}{ll}
\sum_{\sigma\in\Gal(\Qn/\Q)} f(d_n^\sigma) \sigma & \hbox{if $(-1)^n = \varepsilon$} \\
\sum_{\sigma\in\Gal(\Qn/\Q)} f(d_{n-1}^\sigma) \sigma & \hbox{if $(-1)^n = -\varepsilon$}
\end{array}\right. 
$$ 
is an isomorphism from $\Hom(E^\varepsilon(\Qnp),\Zp)$ to $\omega_n^{-\varepsilon} \Zp[\Gal(\Qn/\Q)]$, 
and that for $m \ge n\ge 1$ these maps are compatible in the 
sense that the following diagram commutes
$$
\begin{array}{ccc}
\Hom(E^\pm(\Q_{m,p}),\Zp) &\stackrel{\sim}{\lrar}& \omega_m^\mp \Zp[\Gal(\Q_m/\Q)]\\[4pt]
\big\downarrow&&\big\downarrow\\[4pt]
\Hom(E^\pm(\Qnp),\Zp) &\stackrel{\sim}{\lrar}& \omega_n^\mp \Zp[\Gal(\Qn/\Q)].
\end{array}
$$
Here the left-hand vertical map is restriction, and the right-hand 
vertical map sends $\omega_m^\mp$ to $\omega_n^\mp$.

In the limit it follows (\cite{Ko} Theorem 6.2) that 
$\Hom(E^\pm(\Qinfp),\Zp)$ is free of rank one over $\Iwq$ with a generator 
$f_\pm$ satisfying 
$\sum_{\sigma\in\Gal(\Qn/\Q)} f_\pm(d_n^\sigma) \sigma = \omega_n^\mp$.  
If we take $\mu^\pm$ to be the map corresponding to $f_\pm$ under \eqref{compact}, 
then $\mu^\pm$ satisfies the conclusions of the theorem.
\hfill\qed\vglue12pt

Let $\L_E^\pm \in \Iwq$ denote the $p$-adic $L$-functions defined by the first author in 
Section~6.2.2 of \cite{Po}.  These are characterized by the formulas
\begin{eqnarray}
\label{evenint}
\qquad\; \chi(\L_E^+)= (-1)^{(n+1)/2} 
    \frac{\tau(\chi)}{\chi(\omega_{n}^+)}\frac{L(E,\bar\chi,1)}{\Omega_E}
    &&\hskip-16pt \hbox{if $\chi$ has order $p^n$ with $n$ odd,}\\
\label{oddint}
\chi(\L_E^-)= (-1)^{n/2+1} 
    \frac{\tau(\chi)}{\chi(\omega_{n}^-)}\frac{L(E,\bar\chi,1)}{\Omega_E}
    &&\hskip-16pt\hbox{if $\chi$ has order $p^n>1$ with $n$ even.}
\end{eqnarray}
In addition, if $\chi_0$ is the trivial character then
\begin{equation}
\label{trival}
\chi_0(\L_E^+) = (p-1)\frac{L(E,1)}{\Omega_E}, \quad 
    \chi_0(\L_E^-) = 2\frac{L(E,1)}{\Omega_E}.\hskip.35in
\end{equation}

\proclaim{Theorem}
\label{analytic}
There is an isomorphism $\red{\U}/(\redd{\V}^\pm + \red{\C}) \isom \IwOq/\L_E^\pm\IwOq$.
\endproclaim

{\it Proof}.
By \eqref{lulp} and \eqref{qtok} we have 
\begin{eqnarray*}
\red{\U}/\redd{\V}^\pm &\cong &\Hom_\O(E^\pm(\KQinfp) \otimes \Qp/\Zp, \Kp/\Op) \\
    &\cong &\Hom(E^\pm(\Qinfp) \otimes \Qp/\Zp, \Kp/\Op) \\
    &\cong &\Hom(E^\pm(\Qinfp) \otimes \Qp/\Zp, \Qp/\Zp) \otimes \Op.
\end{eqnarray*}
Let $\mu^\pm$ be as in Theorem \ref{genmu},  
let $\xi$ be the generator of $\rc{\C}$ from Theorem \ref{erl}, and let 
$\varphi^\pm$ be the image of $\xi$ in 
$\Hom_\O(E^\pm(\KQinfp) \otimes \Qp/\Zp, E[p^\infty])$.  
For some $h^\pm  \in \IwOq$ we have 
\begin{equation}
\label{index}
\varphi^\pm = h^\pm \mu^\pm, 
\end{equation}
and then $\red{\U}/(\redd{\V}^\pm + \red{\C})\cong \IwOq/h^\pm\IwOq$.

It follows from \eqref{index} that for every $k, n \ge 1$ and every nontrivial
character $\chi : \Gamma^+ \to \bmu_{p^n}$, 
$$
\sum_{\sigma\in\Gal(\Qn/\Q)} \chi(\sigma) \varphi^\pm(d_n^\sigma \otimes p^{-k}) 
    = \chi(h^\pm)
   \sum_{\sigma\in\Gal(\Qn/\Q)} \chi(\sigma) \mu^\pm(d_n^\sigma \otimes p^{-k}). 
$$
Using the formulas of Theorems \ref{kob}(iv) and \ref{erl} to compute the left-hand side, 
and Theorem \ref{genmu} for the right-hand side, we deduce that if the order 
of $\chi$ is $p^n > 1$ and $\varepsilon = (-1)^{n+1}$ then
$$
\frac{L(E,\bar\chi,1)}{\Omega_E} (-1)^{[\frac{n}{2}]}\tau(\chi) 
    \equiv \chi(h^\varepsilon) \chi(\omega_n^{\varepsilon}) \pmod{p^k}
$$
for every $k$.  It follows from \eqref{evenint} and \eqref{oddint} that 
$h^\pm = -\L_E^\pm$.
\hfill\qed\vglue8pt

The following theorem is our main result.

\proclaim{Theorem}
\label{mainthm}
$\cp_{\Iwq}(\Hom(\Sel^\pm(E/\Qinf),\Qp/\Zp)) = \L_E^\pm \Iwq$.
\endproclaim

{\it Proof}.
We have
\begin{eqnarray*}
\cp_{\IwOq}(\Hom_\O(\Sel^\pm(E/\KQinf),\Kp/\Op)) 
    &= &\cp_{\IwOq}(\red{\X}/\alpha(\redd{\V}^\pm)) \\
    &=& \cp_{\IwOq}(\red{\U}/(\redd{\V}^\pm + \red{\C})) \\
    &= &\L_E^\pm \IwOq
\end{eqnarray*}
by Theorems \ref{newsel}, \ref{bettermc}, and \ref{analytic}, respectively.
Since 
$$
\Sel^\pm(E/\KQinf) = \Sel^\pm(E/\Qinf) \otimes \Op,
$$
we also have
\begin{eqnarray*}
\Hom_\O(\Sel^\pm(E/\KQinf),\Kp/\Op) &=& \Hom(\Sel^\pm(E/\Qinf),\Kp/\Op) \\
    &= &\Hom(\Sel^\pm(E/\Qinf),\Qp/\Zp) \otimes \Op
\end{eqnarray*}
and the theorem follows.
\hfill\qed

\section{Applications}
\label{appl}
\advance\eqcount by 13

We describe briefly the basic applications of the supersingular 
main conjecture.  As in the previous sections, we assume that 
$E$ is an elliptic curve defined over $\Q$, with complex multiplication 
by the ring of integers of an imaginary quadratic field $K$, and $p$ 
is an odd prime where $E$ has good supersingular reduction.  
For this section we write $\Gamma = \Gamma^+$, so $\Lambda = \Zp[[\Gamma]]$.

\numbereddemo{{R}emark}
The results below also hold for primes of ordinary reduction, 
and can be proved using the main conjecture for ordinary primes.
\enddemo 

The following application was already proved in \cite{Ru3}, 
as an application of Theorem \ref{oldmc}.

\proclaimtitle{[Ru3, Th.~11.4]}
\proclaim{Theorem}
\label{app1}
If $L(E,1) \ne 0${\rm ,} then $E(\Q)$ is finite and 
$$
|\Sh(E)| = r \frac{L(E,1)}{\Omega_E}
$$
where $r \in \Q^\times$ satisfies $\ord_p(r) = 0${\rm ,} 
as predicted by the Birch and Swinnerton\/{\rm -}\/Dyer conjecture.

If $L(E,1) = 0${\rm ,} then either $E(\Q)$ is infinite or 
$\Sh(E)[p^\infty]$ is infinite.
\endproclaim

Before proving Theorem \ref{app1} we need the following lemma.

\proclaim{Lemma}
\label{cont}
The natural restriction map 
$\Sel(E/\Q) \to \Sel^\pm(E/\Qinf)^\Gamma$ is an isomorphism.
\endproclaim

{\it Proof}.
For every number field $F$ let $\Sel'(E/F)$ denote the Selmer group of $E$ 
over $F$ with no local condition at $p$:
$$
\Sel'(E/F) = \ker : H^1(F,E[p^\infty]) \to 
    \oplus_{v \nmid p} H^1(F_v,E[p^\infty])
$$
(note that $E(F_v) \otimes \Qp/\Zp = 0$ when $v \nmid p$).  Thus 
we have a commutative diagram
\begin{equation}
\label{snake}
\begin{array}{cccccccc}
0 & \too & \Sel(E/\Q) & \too & \Sel'(E/\Q) & \too & H^1(\Qp,E[p^\infty])/A \\
& & \downarrow  && \downarrow  && \downarrow \\
0 & \too & \Sel^\pm(E/\Qinf)^\Gamma & \too & \Sel'(E/\Qinf)^\Gamma & \too & 
 H^1(\Qinfp,E[p^\infty])/A^\pm_\infty
\end{array}
\end{equation}
where $A$ and $A^\pm_\infty$ are the images of $E(\Qp) \otimes \Qp/\Zp$ and 
$E^\pm(\Qinfp) \otimes \Qp/\Zp$, respectively, and the vertical maps are 
restriction maps.
It follows from the theory of complex multiplication that 
$E(\Qinfp)$ has no $p$-torsion, and then standard methods (see for example 
Proposition 1.2 of \cite{Ru1}) 
show that the restriction maps 
$$
H^1(\Qp,E[p^\infty]) \to H^1(\Qinfp,E[p^\infty])^\Gamma, \quad 
\Sel'(E/\Q) \to \Sel'(E/\Qinf)^\Gamma
$$ 
are isomorphisms.

We will show that for every $n$ the map 
$
E(\Qp)\otimes \Qp/\Zp \to (E^\pm(\Qnp)\otimes \Qp/\Zp)^\Gamma
$ 
is surjective.  It will then follow that 
the right-hand vertical map in \eqref{snake} is injective, and then 
(using the remarks above and the snake lemma) that the left-hand vertical map 
in \eqref{snake} is an isomorphism, which is the assertion of the lemma.

To show that 
$
E(\Qp)\otimes \Qp/\Zp \to (E^\pm(\Qnp)\otimes \Qp/\Zp)^\Gamma
$ 
is surjective  
it suffices to check that $\dim_{\Fp}(E^\pm(\Qnp)\otimes\Fp)^\Gamma = 1$, 
since $E(\Qp)\otimes \Qp/\Zp \cong \Qp/\Zp$.
Identify $\Fp[\Gal(\Qn/\Q)]$ with $\Fp[X]/(X^{p^n}-1) = \Fp[X]/(X-1)^{p^n}$. 
Since $E^\pm(\Qnp)$ is cyclic over $\Zp[\Gal(\Qn/\Q)]$ (Theorem \ref{kob}(v)), 
$$
E^\pm(\Qnp)\otimes\Fp \cong \Fp[X]/(X-1)^a 
$$
for some $a \ge 0$.  
Under this identification $(E^\pm(\Qnp)\otimes\Fp)^\Gamma$ is the kernel of 
$X-1$, which is visibly \pagebreak one-dimensional.
\hfill\qed

\demo{{P}roof of Theorem {\rm 8.2}}
By Lemma \ref{cont} we have
\begin{eqnarray*}
|\Sel(E/\Q)| &=& |\Hom(\Sel(E/\Q),\Qp/\Zp)| \\
    &=& |\Hom(\Sel^\pm(E/\Qinf)^\Gamma,\Qp/\Zp)|\\
    &=& |\Hom(\Sel^\pm(E/\Qinf),\Qp/\Zp) \otimes_\Lambda \Zp|.
\end{eqnarray*}
By Theorems \ref{newsel} and \ref{bettermc}, $\Hom(\Sel^\pm(E/\Qinf),\Qp/\Zp)$ 
has no nonzero finite submodules, and by Theorem \ref{mainthm} its 
characteristic ideal is $\L_E^\pm\Lambda$.  Writing $\chi_0$ for the trivial character 
of $\Gamma$, standard techniques (for example \cite[Lemma~4 of \S{I.1.3}]{PR1}) 
show that 
$$
|\Hom(\Sel^\pm(E/\Qinf),\Qp/\Zp) \otimes_\Lambda \Zp| = |\Zp/\chi_0(\L_E^\pm)\Zp| 
    = |\Zp/(L(E,1)/\Omega_E)\Zp|
$$
using \eqref{trival} for the last equality.  This proves the theorem.
\hfill\qed\enddemo

Fix a generator $\gamma$ of $\Gamma$.  Define $\nu_0 = \gamma-1$ 
and for every $n \ge 1$ let $\nu_n = \sum_{i=0}^{p-1} \gamma^{ip^{n-1}}$.  
If $\chi$ is a character of $\Gamma$ of finite order, 
let $\Zp[\chi]$ denote the ring obtained by adjoining the values 
of $\chi$ to $\Zp$.  We view $\Zp[\chi]$ as a $\Lambda$-module with $\Gamma$ 
acting via $\chi$, and if $M$ is a $\Lambda$-module
we define $M^\chi = M \otimes_{\Lambda} \Zp[\chi]$.
Then $\chi(\nu_m) = 0$ if and only if the order of $\chi$ is $p^m$, and 
if $M$ is finitely generated or co-finitely generated over $\Zp$ 
and $\chi$ has order $p^m$, then 
$M^\chi$ is infinite if and only if $M^{\nu_m = 0}$ is infinite, 
where $M^{\nu_m = 0}$ is the kernel of $\nu_m$ on $M$.

For every $n$ write $G_n = \Gal(\Qn/\Q)$.  

\proclaim{Theorem}
\label{app2}
Suppose $\chi$ is a character of $G_n$.  
If $L(E,\chi,1) \ne 0$ then $E(\Qn)^\chi$ and $\Sh(E/\Qn)^\chi$ are finite.  
If $L(E,\chi,1) = 0$ then either
$E(\Qn)^\chi$ is infinite or $\Sh(E/\Qn)^\chi$ is infinite.
\endproclaim

Before proving Theorem \ref{app2} we need the following lemma. 

\proclaim{Lemma}
\label{lema2}
Suppose $\chi$ is a character of $G_n$ of order $p^m > 1${\rm ,} and let 
$\varepsilon = (-1)^m$.  Then $\Sel^\varepsilon(E/\Qn)^{\nu_m = 0}$ is infinite if and 
only if $\Sel(E/\Qn)^{\nu_m = 0}$ is infinite. 
\endproclaim

{\it Proof}.
We have $\Sel^\varepsilon(E/\Qn) \subset \Sel(E/\Qn)$, so one implication is clear.  
Suppose now that $\Sel(E/\Qn)^{\nu_m = 0}$ is infinite.  By  
Proposition 10.1 of \cite{Ko}, either $\Sel^\varepsilon(E/\Qn)^{\nu_m = 0}$ 
or $\Sel^{-\varepsilon}(E/\Qn)^{\nu_m = 0}$ must be infinite.  
But localization at $p$ sends $\Sel^{-\varepsilon}(E/\Qn)^{\nu_m = 0}$ into 
$E^{-\varepsilon}(\Qnp)^{\nu_m = 0}$ which is zero, and so 
$\Sel^{-\varepsilon}(E/\Qn)^{\nu_m = 0} \subset \Sel^\varepsilon(E/\Qn)^{\nu_m = 0}$.
Hence $\Sel^\varepsilon(E/\Qn)^{\nu_m = 0}$ is infinite.
\hfill\qed

\demo{Proof of Theorem {\rm \ref{app2}}}
Let $p^m$ be the order of $\chi$.  If $m = 0$ then the theorem is a 
consequence of Theorem \ref{app1}.  So we may suppose $m \ge 1$, and 
we let $\varepsilon = (-1)^m$.
\begin{eqnarray*}
\hbox{$\Sel(E/\Qn)^\chi$ is infinite} &\iff &
    \hbox{$\Sel(E/\Qn)^{\nu_m = 0}$ is infinite} \\ &\iff &
    \hbox{$\Sel^\varepsilon(E/\Qn)^{\nu_m = 0}$ is infinite} \\ &\iff &
    \hbox{$\Sel^\varepsilon(E/\Qinf)^{\nu_m = 0}$ is infinite} \\ &\iff &
    \hbox{$\Hom(\Sel^\varepsilon(E/\Qinf),\Qp/\Zp) \otimes \Lambda/\nu_m$ is infinite} \\ &\iff &
    \hbox{$\Lambda/(\L^\varepsilon_E,\nu_m)$ is infinite} \\ &\iff&
    \bar\chi(\L^\varepsilon_E) = 0 \\ &\iff&
    L(E,\chi,1) = 0
\end{eqnarray*}
using Lemma \ref{lema2}, Theorem 9.3 of \cite{Ko}, 
Theorem \ref{mainthm}, and \eqref{evenint} and \eqref{oddint}.
\enddemo
 
\vglue-12pt
 \AuthorRefNames [999]

\end{document}